\documentclass[a4paper, 11pt]{article}
\usepackage{amsmath} \usepackage{amsfonts} \usepackage{amssymb}\usepackage{latexsym}
\usepackage[english]{babel}
\usepackage{amscd}
\usepackage[dvips,final]{graphics}
\usepackage{locengli,math}
\usepackage{graphicx,pst-plot,psfig,pstricks,pst-node,pst-text,pst-tree,pst-char,pst-coil}
\begin{document}
\begin{center}
\textbf{\LARGE{\textsf{Construction of Nijenhuis operators and dendriform trialgebras}}}
\vskip1cm
\parbox[t]{14cm}{\large{
Philippe {\sc Leroux}}\\
\vskip4mm
{\footnotesize
\baselineskip=5mm
Institut de Recherche
Math\'ematique, Universit\'e de Rennes I and CNRS UMR 6625\\
Campus de Beaulieu, 35042 Rennes Cedex, France, pleroux@univ-rennes1.fr}}
\end{center}
\vskip1cm
\baselineskip=5mm
\noindent
\begin{center}
\begin{large}
\textbf{ 05/11/03}
\end{large}
\end{center}
{\bf Abstract:}
We construct Nijenhuis operators from particular bialgebras called dendriform-Nijenhuis bialgebras. It turns out that such Nijenhuis operators commute with $TD$-operators, kind of Baxter-Rota operators, and therefore closely related
dendriform trialgebras. This allows the construction of associative algebras, called dendriform-Nijenhuis algebras made out with nine operations and presenting an exotic combinatorial property. We also show that the augmented free dendriform-Nijenhuis algebra and its commutative version have a structure of connected Hopf algebras. Examples are given.
\section{Introduction}
\textbf{Notation:}
In the sequel $k$ is a field of characteristic zero.
Let $(X, \ \diamond)$ be a $k$-algebra and $(\diamond_i)_{1 \leq i \leq N}: X^{\otimes 2} \xrightarrow{} X$ be a family of binary operations on $X$. 
The notation $\diamond \longrightarrow \sum_i \diamond_i$ will mean $x \diamond y =  \sum_i x \diamond_i y$, for all $x,y \in X$. We say that
the operation $\diamond$ {\it{splits}} into the $N$ operations $\diamond_1, \ldots, \diamond_N,$ or that the operation $\diamond$ is a {\it{cluster of $N$ (binary) operations}}. 

Let $(\mathcal{L}, \ [,])$ be a Lie algebra. A {\it{Nijenhuis operator}} $N:\mathcal{L} \xrightarrow{} \mathcal{L}$ is a linear map verifying,
$$ [N(x),N(y)] + N^2([x,y]) = N([N(x),y] + [x, N(y)]).    \ \ \ (1)$$
Solutions of this equation can be constructed by producing operators $\beta: A \xrightarrow{} A$ verifying the {\it{associative Nijenhuis relation (ANR)}} on an associative algebra $(A, \ \mu)$. 
The {\it{associative Nijenhuis relation}}, i.e., 
$$(ANR): \ \ \ \ \beta(x)\beta(y) + \beta^2(xy) = \beta(\beta(x)y + x\beta(y)), \ \ \ \ (2)$$
for all $x,y \in A$, appears for the first time in \cite{CGG}, see also  \cite{KEF, KEFNij} and the references therein. Such linear maps $\beta$ are then Nijenhuis operators since $(2)$ implies $(1)$ on the Lie algebra $(A, \ [,])$, with $[x,y]= xy-yx$, for all $x,y \in A$. In the sequel, by Nijenhuis operators, we mean a linear map defined on an associative algebra and verifying $(2)$.

Section~2 prepares the sequel of this work. $NS$-Algebras are defined and dendriform trialgebras are recalled. We show that Nijenhuis operators on an associative algebra give $NS$-algebras. The notion of $TD$-operators is also introduced. Such operators give dendriform trialgebras \cite{LodayRonco, KEF, Lertribax}. In Section~3, the notion of dendriform-Nijenhuis bialgebras
$(A, \ \mu, \ \Delta)$ is introduced. This notion is the corner-stone of the paper. Indeed, any dendriform-Nijenhuis bialgebra gives two operators $\beta, \ \gamma: \textsf{End}(A) \xrightarrow{} \textsf{End}(A)$, where $\textsf{End}(A)$ is the $k$-algebra of linear maps from $A$ to $A$, which commute one another, i.e., $\beta\gamma = \gamma\beta$. The first one turns out to be a Nijenhuis operator and the last one, a $TD$-operator. Section~4 gives examples. Section~5 introduces dendriform-Nijenhuis algebras, which are associative algebras whose associative product splits into nine binary operations linked together by 28 constraints. When the nine operations are gathered in a particular way, dendriform-Nijenhuis algebras become $NS$-algebras and when gathered in another way, dendriform-Nijenhuis algebras become
dendriform trialgebras. Otherwise stated, this means the existence of a commutative diagram between the involved categories,
\begin{center}

$
\begin{array}{ccc}
\textsf{Dend. Nijenhuis} & \longrightarrow & \textsf{NS} \\ 
 \downarrow & \searrow     & \downarrow    \\
\textsf{TriDend.} & \longrightarrow  & \textsf{As.} 
\end{array} 
$
\end{center}
diagram which will be explained in Section~5.
Such associative algebras can be easily constructed from dendriform-Nijenhuis bialgebras. Section~6 shows that the augmented free dendriform-Nijenhuis algebra and the augmented free $NS$-algebra on a $k$-vector space $V$, as well as their commutative versions, have a structure of connected Hopf algebras.
\section{$NS$-Algebras and dendriform trialgebras}
We present some
results relating $NS$-algebras and dendriform trialgebras to Nijenhuis operators and $TD$-operators. 
\subsection{$NS$-Algebras}
\begin{defi}{[$NS$-algebra]}
A {\it{$NS$-algebra}} $(A, \ \prec, \ \succ, \ \bullet )$ is a $k$-vector space equipped with three binary operations $\prec, \ \succ, \ \bullet: A^{\otimes 2} \xrightarrow{} A$ verifying,
\begin{eqnarray*}
(x \prec y) \prec z  := x \prec ( y \star z ), \ \ 
(x \succ y) \prec z  &:=& x \succ (y \prec z ) , \ \ 
(x \star y) \succ z := x \succ y \succ z, \\
(x\star y)\bullet z + (x \bullet y) \prec z &:= & x \succ (y \bullet z) + x \bullet (y \star z),
\end{eqnarray*}
where $\star \longrightarrow \ \prec + \succ + \ \bullet$. The $k$-vector space $(A, \ \star)$ is an associative algebra, i.e., there exits a functor $F_1: \textsf{NS} \longrightarrow \textsf{As.}$, where $\textsf{NS}$ is the category of $NS$-algebras and $\textsf{As.}$ the category of associative algebras.
\end{defi}
\Rk
Let  $(A, \ \prec, \ \succ, \ \bullet )$ be a $NS$-algebra. Define three operations $\prec^{op}, \ \succ^{op}, \ \bullet^{op}: A^{\otimes 2} \xrightarrow{} A$ as follows,
$$ x \prec^{op} y := y \succ x, \ \ x \succ^{op} y := y \prec x, \ \ x \bullet^{op} y := y \bullet x, \ \ \forall x,y \in A.$$
Then, $(A, \ \prec^{op}, \ \succ^{op}, \ \bullet^{op})$ is a $NS$-algebra called the {\it{opposite}} of $(A, \ \prec, \ \succ, \ \bullet)$. $NS$-Algebras are said to be {\it{commutative}} if they coincide with their opposite, i.e., if $x \prec y = y \succ x$ and  $x \bullet y := y \bullet x$.
\begin{prop}
\label{asnij}
Let $(A, \ \mu)$ be an associative algebra equipped with a Nijenhuis operator $\beta: A \xrightarrow{} A$. Define three binary operations $\prec_{\beta}, \ \succ_{\beta}, \ \bullet_{\beta}: A^{\otimes 2} \xrightarrow{} A$ as follows,
$$ x \prec_{\beta} y := x \beta(y), \ \  x \succ_{\beta} y := \beta(x) y, \ \  x \bullet_{\beta} y := -\beta(xy), $$ 
for all  $x,y \in A$.
Then,
$A^{\beta} := (A, \ \prec_{\beta}, \ \succ_{\beta}, \ \bullet_{\beta} )$ is a $NS$-algebra.  
\end{prop}
\Proof
Straightforward.
\eproof
\begin{prop}
\label{gamnin}
Let $(A, \ \mu)$ be a unital associative algebra with unit $i$. Suppose  $\beta: A \xrightarrow{} A$ is a Nijenhuis operator. Define three binary operations $\tilde{\prec}_{\beta}, \ \tilde{\succ}_{\beta}, \ \tilde{\bullet}_{\beta}: (A^\beta)^{\otimes 2} \xrightarrow{} A^\beta$ as follows,
$$ x \tilde{\prec}_{\beta} y := x \beta(i)\prec_{\beta} y = x \prec_{\beta} \beta(i)y, \ \  x \tilde{\succ}_{\beta} y :=  x \succ_{\beta} \beta(i) y = x \beta(i) \succ_{\beta}  y, \ \  x \tilde{\bullet}_{\beta} y := x \beta(i)\bullet_{\beta} y = x \bullet_{\beta} \beta(i) y,$$ 
for all $x,y \in A.$ 
Then,
$(A, \ \tilde{\prec}_{\beta}, \ \tilde{\succ}_{\beta}, \ \tilde{\bullet}_{\beta} )$ is a $NS$-algebra. 
\end{prop}
\Proof
Observe that for all $x \in A$, $\beta(x) \beta(i) = \beta(x \beta(i))$ and $\beta(i)\beta(x)  = \beta(\beta(i) x )$. Therefore, for all $x,y \in A$, $x \tilde{\prec}_{\beta} y := x \beta(i)\beta(y) = x \beta(\beta(i)y) = x \prec_{\beta} \beta(i)y $ and similarly for the two other operations. Fix $x,y,z \in A$. Let us check that,
$$(x \tilde{\prec}_{\beta} y ) \tilde{\prec} z = x \tilde{\prec}_{\beta} (y  \tilde{\star}_{\beta} z),$$
where, $\tilde{\star}_{\beta} \longrightarrow \ \tilde{\prec}_{\beta} + \tilde{\succ}_{\beta} + \ \tilde{\bullet}_{\beta}$.
Indeed,
\begin{eqnarray*}
(x \tilde{\prec}_{\beta} y ) \tilde{\prec}_{\beta} z &:=& x \prec_{\beta} (\beta(i) y  \star_{\beta} \beta(i) z), \\
&=& x \prec_{\beta} (\beta(i) y  \prec_{\beta} \beta(i) z + \beta(i) y  \succ_{\beta} \beta(i) z + \beta(i) y  \bullet_{\beta} \beta(i) z), \\
&=& x \prec_{\beta} \beta(i) (y  \prec_{\beta} \beta(i) z +  y  \succ_{\beta} \beta(i) z +  y  \bullet_{\beta} \beta(i) z), \\
&=& x \tilde{\prec}_{\beta} (y  \tilde{\star}_{\beta} z).
\end{eqnarray*}
\eproof
\subsection{Dendriform trialgebras}
Let us recall some motivations for the introduction of dendriform trialgebras.
Motivated by $K$-theory, J.-L. Loday first introduced  a ``non-commutative version'' of Lie algebras called Leibniz algebras \cite{LodayLeib}. Such algebras are described by a bracket $[ -,z]$ verifying the Leibniz identity:
$$[[x,y ] ,z ] = [[ x,z ] ,y ] + [ x,[ y  ,z ]].$$
When the bracket is skew-symmetric, the Leibniz identity becomes the
Jacobi identity and Leibniz algebras turn out to be Lie algebras. A way to construct such  Leibniz algebras is to start with  associative dialgebras, i.e., $k$-vector spaces $D$ equipped with two associative products,
$\vdash$ and $\dashv$, and verifying some conditions \cite{Loday}. The operad $Dias$ associated with associative dialgebras is then Koszul dual to the operad $DiDend$ associated with dendriform dialgebras \cite{Loday}.
A {\it{dendriform dialgebra}} is a $k$-vector space $E$ equipped with two binary operations:
$\prec, \ \succ: E^{\otimes 2} \xrightarrow{} E$, satisfying the following relations for all $x,y \in E$:
$$(x \prec y )\prec z = x \prec(y \star z), \ \ \
(x \succ y )\prec z = x \succ(y \prec z), \ \ \
(x \star y )\succ z = x \succ(y \succ z), \ \ \ $$
where by definition $x \star y :=x  \prec y +x \succ y$, for all $x,y \in E$. The dendriform dialgebra $(E,\star)$ is then an associative algebra such that $\star \longrightarrow \ \prec + \succ$. 
Similarly, to propose a ``non-commutative version'' of Poisson algebras, J.-L. Loday and M. Ronco \cite{LodayRonco} introduced the notion of associative trialgebras. It turns out that $Trias$, the operad associated with this type of algebras, is Koszul dual to $TriDend$, the operad associated with
dendriform trialgebras.
\begin{defi}{[Dendriform trialgebra]}
\label{tridend}
A {\it{dendriform trialgebra}} is a $k$-vector space $T$ equipped with three binary operations:
$\prec, \ \succ, \ \circ: T^{\otimes 2} \xrightarrow{} T$, satisfying the following relations for all $x,y \in T$:
$$(x \prec y )\prec z = x \prec(y \star z), \ 
(x \succ y )\prec z = x \succ(y \prec z), \
(x \star y )\succ z = x \succ(y \succ z), $$
$$(x \succ y )\circ z = x \succ(y \circ z), \ 
(x  \prec y )\circ z = x \circ(y \succ z), \ 
(x \circ y )\prec z = x \circ(y \prec z), \
(x \circ y ) \circ z = x \circ(y \circ z),$$
where by definition $x \star y :=x  \prec y +x \succ y +x \circ y$, for all $x,y \in T$. The $k$-vector space $(T,\star)$ is then an associative algebra such that $\star \longrightarrow \ \prec + \succ + \ \circ$. There exists a functor $F_2: \textsf{TriDend.} \xrightarrow{} \textsf{As}.$
\end{defi}
Observe that these axioms are globally invariant under the transformation $x \prec^{op} y := y \succ x$, $x \succ^{op} y := y \prec x$ and $x \circ^{op} y := y \circ x$. A dendriform trialgebra is said to be {\it{commutative}} if $x \prec y := y \succ x$ and $x \circ y := y \circ x$.

\noindent
To construct dendriform trialgebras, the use of $t$-Baxter operators, also called Rota-Baxter operators have been used in \cite{KEF} and in \cite{Lertribax} to generalise \cite{AguiarLoday}.
Let $(A, \ \mu)$ be an associative algebra and $t \in k$. A {\it{$t$-Baxter operator}} is a linear map $\xi: A \xrightarrow{} A$ verifying:
$$\xi(x) \xi(y) = \xi(x  \xi(y) + \xi(x) y + txy).  $$
For $t=0$, this map is called a Baxter operator. It 
appears originally in a work of G. Baxter \cite{Baxter} and the importance of such a map was stressed by G.-C. Rota in \cite{Rota1}.
We present another way to produce dendriform trialgebras.
Let $(A, \ \mu)$ be a unital associative algebra with unit $i$. The linear map $\gamma: A \xrightarrow{} A$ is said to be a {\it{$TD$-operator}} if,
$$\gamma(x)\gamma(y) = \gamma(\gamma(x)y + x\gamma(y) -x \gamma(i)y), $$
for all $x,y \in A$.
\begin{prop}
\label{TDdend}
Let $A$ be a unital algebra with unit $i$. Suppose $\gamma: A \xrightarrow{} A$ is a $TD$-operator. Define three binary operations $\prec_{\gamma}, \ \succ_{\gamma}, \ \circ_{\gamma}: A^{\otimes 2} \xrightarrow{} A$ as follows,
$$ x \prec_{\gamma} y := x \gamma(y), \ \  x \succ_{\gamma} y := \gamma(x) y, \ \  x \circ_{\gamma} y := -x \gamma(i) y, \ \ \forall \ x,y \in A.$$ Then,
$A^{\gamma}:=(A, \ \prec_{\gamma}, \ \succ_{\gamma}, \ \circ_{\gamma} )$ is a dendriform trialgebra. The operation $\bar{\star}_{\gamma}: A^{\otimes 2} \xrightarrow{} A$ defined by,
$$ x \bar{\star}_{\gamma} y := x \gamma(y) +  \gamma(x) y - x \gamma(i) y,$$
is associative.
\end{prop}
\Proof
Straightforward by noticing that $\gamma(i)\gamma(x) = \gamma(x)\gamma(i)$, for all $x \in A$.
\eproof
\section{Construction of Nijenhuis operators and $TD$-operators from dendriform-Nijenhuis bialgebras}
In \cite{Aguiar}, Baxter operators are constructed from $\epsilon$-bialgebras. This idea have been used in 
 \cite{AguiarLoday, Lertribax} to produce commuting $t$-Baxter operators.
Recall that a {\it{$t$-infinitesimal bialgebra}} (abbreviated $\epsilon(t)$-bialgebra) is a triple $(A, \ \mu, \ \Delta )$ where $(A, \ \mu)$ is an associative algebra and $(A,\ \Delta )$ is a coassociative coalgebra such that for all $a,b \in A$,
$$\Delta (ab)= a_{(1)} \otimes a_{(2)} b + a b_{(1)} \otimes b_{(2)}  + t a \otimes b.$$
If $t=0$, a $t$-infinitesimal bialgebra is called
an {\it{infinitesimal bialgebra}} or a $\epsilon$-bialgebra. Such  bialgebras
appeared for the first time in the work of Joni and Rota in \cite{Rota}, see also Aguiar \cite{Aguiar}, for the case $t=0$ and Loday \cite{Lodayscd}, for the case $t=-1$.

\noindent
To produce Nijenhuis operators from bialgebras, we replace the term $t a \otimes b$ in the definition of $t$-infinitesimal bialgebras by
the term $-\mu(\Delta(a)) \otimes b$. A physical interpretation of this term can be the following. The two-body system $\mu(a \otimes b) := ab$, for instance a two particule system or a string of letters in informatic, made out from $a$ and $b$ is sounded by a coproduct $\Delta$, representing a physical system. This coproduct ``reads sequentially" the system $ \mu(a \otimes b)$ giving the systems $\Delta(a)b$ and $a \Delta(b)$. In the case of $t$-infinitesimal bialgebras, $\Delta(ab)$ studies the behavior between the "sequential reading", $\Delta(a)b$ and $a \Delta(b)$, and $t a \otimes b$ which can be interpreted as the system $a$ decorrelated with system $b$. In the replacement, $t a \otimes b$ by $\mu(\Delta(a)) \otimes b$, we want to compare the "sequential reading", $\Delta(a)b$ and $a \Delta(b)$ to the decorrelated system $\mu(\Delta(a)) \otimes b$ made out with $b$ and the system obtained from the recombinaison of the pieces created by reading the system $a$.
\begin{defi}{[Dendriform-Nijenhuis bialgebra]}
A {\it{Dendriform-Nijenhuis bialgebra}} is a triple $(A, \ \mu, \ \Delta )$ where $(A, \ \mu)$ is an associative algebra and $(A,\ \Delta )$ is a coassociative coalgebra such that,
$$ \Delta(ab) := \Delta(a)b + a\Delta(b) - \mu(\Delta(a)) \otimes b, \ \forall a,b \in A.$$
\end{defi}
The $k$-vector space $\textsf{End}(A)$ of linear
endomorphisms of $A$ is viewed as an associative algebra under composition denoted simply by concatenation $TS$, for $T,S \in \textsf{End}(A)$. Another operation called the convolution product, $*$, defined by $T *S := \mu(T \otimes S)\Delta$, for all $T,S \in \textsf{End}(A)$ will be used.

\begin{prop}
\label{NTDop}
Let $(A, \ \mu, \ \Delta)$ be a dendriform-Nijenhuis bialgebra. Equip $\textsf{End}(A)$ with the convolution product $*$. Define the operators $\beta, \gamma: \textsf{End}(A) \xrightarrow{} \textsf{End}(A)$ by $T \mapsto \beta(T):= id * T$ (right shift) and $T \mapsto \gamma(T):= T*id $ (left shift). Then, the right shift
$\beta$ is a Nijenhuis operator and the left shift $\gamma$ is a $TD$-operator. Moreover $\beta\gamma=\gamma\beta$.
\end{prop}
\Proof
Let us show that the right shift $\beta: T \mapsto id * T$ is a Nijenhuis operator.
Fix $T, S \in \textsf{End}(A)$ and $a \in A$. On the one hand,
\begin{eqnarray*}
\beta(T)\beta(S)(a) &=& \mu(id \otimes T) \Delta (a_{(1)} S(a_{(2)})), \\
&=& \mu(id \otimes T) (a_{(1)(1)} \otimes a_{(1)(2)} S(a_{(2)}) + a_{(1)}
S(a_{(2)})_{(1)}  \otimes S(a_{(2)})_{(2)} -  a_{(1)(1)} a_{(1)(2)} \otimes  S(a_{(2)})),\\
&=& a_{(1)(1)}T( a_{(1)(2)} S(a_{(2)})) + a_{(1)}
S(a_{(2)})_{(1)} T( S(a_{(2)})_{(2)}) -  a_{(1)(1)} a_{(1)(2)} T( S(a_{(2)})).
\end{eqnarray*}
On the other hand,
\begin{eqnarray*}
\beta(\beta(T)S)(a) &=& \mu(id \otimes \beta(T)S) (a_{(1)} \otimes a_{(2)}), \\
&=& a_{(1)}
S(a_{(2)})_{(1)} T( S(a_{(2)})_{(2)})
\end{eqnarray*}
And,
\begin{eqnarray*}
\beta(T \beta(S))(a) &=& \mu(id \otimes T\beta(S)) (a_{(1)} \otimes a_{(2)}), \\
&=& a_{(1)}T( a_{(2)(1)} S(a_{(2)(2)})), \\
&=& a_{(1)(1)}T( a_{(1)(2)} S(a_{(2)})),
\end{eqnarray*}
since $\Delta$ is coassociative. Moreover,
\begin{eqnarray*}
\beta(\beta(TS))(a) &=& \mu(id \otimes \beta(TS)) (a_{(1)} \otimes a_{(2)}), \\
&=& a_{(1)} \beta(TS)(a_{(2)}), \\
&=& a_{(1)} \mu(id \otimes TS) ( a_{(2)(1)} \otimes a_{(2)(2)}), \\
&=& a_{(1)}  a_{(2)(1)} TS(a_{(2)(2)}), \\
&=& a_{(1)(1)}  a_{(1)(2)} TS(a_{(2)}),
\end{eqnarray*}
since $\Delta$ is coassociative. Similarly, we can show that the left shift $\gamma: T \mapsto T * id$ is
a $TD$-operator. 
\eproof

\noindent
A
{\it{$L$-anti-dipterous}} algebra $(A, \ \bowtie, \prec_{A})$ is an associative algebra $(A, \ \bowtie)$ equipped with
a right module on itself i.e., $(x \prec_{A} y) \prec_{A} z=x \prec_{A} (y \bowtie z)$ and such that $\bowtie$ and $\prec_{A}$ are linked by the relation
$(x \bowtie y)\prec_{A}  z = x \bowtie (y \prec_{A} z).$
This notion has been introduced in \cite{BaxLer} and comes from a particular notion of bialgebras. See also \cite{LR1}.
\begin{prop}
Let $(A, \ \mu, \ \Delta)$ be a dendriform-Nijenhuis bialgebra. 
Set $x \bowtie y := \mu(\Delta(x)) y$ and $x \prec_{A} y := x \mu(\Delta(y))$ defined for all $x,y \in A$.
Then, the $k$-vector space $(A, \ \bowtie, \ \prec_{A} )$ is a $L$-anti-dipterous algebra.
\end{prop}
\Proof
Straightforward by using $\mu(\Delta(xy))= x\mu(\Delta(y))$, for all $x,y \in A$.
\eproof

\noindent
Let $(A, \ \mu, \ \Delta)$ be a dendriform-Nijenhuis bialgebra.
We end this subsection by showing that
$A$ can admit a left counit and that 
the $k$-vector space $\textsf{Der}(A)$ of derivatives from $A$ to $A$ is stable by the right shift $\beta$.
\begin{prop}
Let $(A, \ \mu, \ \Delta)$ be a dendriform-Nijenhuis bialgebra. As a $k$-algebra, set $A:=k \bra S \ket  / \mathcal{R}$, where $S$ is a non-empty set and $\mathcal{R}$ is a set of relations. Suppose $\eta: A \xrightarrow{} k$ is a $\mu$-homomorphism and $(\eta \otimes id)\Delta := id$ on $S$. Then, $\eta$ is a left counit, i.e., $(\eta \otimes id)\Delta := id$ on $A$.
\end{prop}
\Proof
Fix $a,b \in A$. Suppose $\eta: A \xrightarrow{} k$ is a $\mu$-homomorphism and $(\eta \otimes id)\Delta (a) := a$ and $(\eta \otimes id)\Delta (b) := b$. Then, 
$$(\eta \otimes id)\Delta (ab) := \eta(a_{(1)})a_{(2)}b + \eta(a)\eta(b_1)b_2 - \eta(a_{(1)})\eta(a_{(2)})b.$$
However, $\eta(a_{(1)})a_{(2)}b := ab$, $\eta(a)\eta(b_1)b_2 := \eta(a)b$ and $\eta(a_{(1)})\eta(a_{(2)})b := \eta(\eta(a_{(1)})a_{(2)})b := \eta(a)b$. Therefore, $(\eta \otimes id)\Delta (ab) := ab$.
\eproof
\begin{prop}
Let $(A, \ \mu, \ \Delta)$ be a dendriform-Nijenhuis bialgebra. If $\partial: A \xrightarrow{} k$ is a derivative, i.e., $\partial(xy):=\partial(x)y + x\partial(y)$, then, so is the linear map $\beta(\partial):x \mapsto \mu(id \otimes \partial)\Delta (x)$. 
\end{prop}
\Proof
Straightforward.
\eproof
\section{Examples}
\begin{prop}
\label{TDN Iff}
Let $S$ be a set. Suppose $\Delta: kS \xrightarrow{} kS^{\otimes 2}$ is a coassociative coproduct on the free $k$-vector space spanned by $S$. Denote by $As(S)$ the free associative algebra generated by $S$ and extend the coproduct $\Delta$ to $\Delta_\sharp: As(S) \xrightarrow{} As(S)^{\otimes 2}$ as follows,
$$  \Delta_\sharp(s) := \Delta(s), \ \ \forall \ s \in kS,$$
$$ \Delta_\sharp(ab) := \Delta_\sharp(a)b + a \Delta_\sharp(b) - \mu(\Delta_\sharp(a))  \otimes b,$$
for all $a,b \in As(S)$.
Then, $(As(S), \ \Delta_\sharp)$ is a dendriform-Nijenhuis bialgebra.
\end{prop}
\Proof
Keep notation of Proposition \ref{TDN Iff}.
The co-operation $\Delta_\sharp$ is well defined since it does not depend on the writting of a given element $c \in As(S)$. Indeed, it is straightforward to show that if $c=ab =a'b' \in As(S)$, with $a,b,a',b' \in As(S)$, then, $$\Delta_\sharp(c):= \Delta_\sharp(a)b + a \Delta_\sharp(b) - \mu(\Delta_\sharp(a))  \otimes b=\Delta_\sharp(a')b' + a' \Delta_\sharp(b') - \mu(\Delta_\sharp(a'))  \otimes b'.$$
Let us show that $\Delta_\sharp$ is coassociative. Let $a,b \in As(S)$. Write $\Delta_\sharp(a) := a_{(1)} \otimes a_{(2)}$ and $\Delta_\sharp(b) := b_{(1)} \otimes b_{(2)}$. Suppose $(id \otimes \Delta_\sharp)\Delta_\sharp(x) = (\Delta_\sharp \otimes id)\Delta_\sharp(x)$, for $x=a, \ b$. By definition,
$\Delta_\sharp(ab) := \Delta_\sharp(a)b + a \Delta_\sharp(b) - \mu(\Delta_\sharp(a))  \otimes b$. On the one hand,
\begin{eqnarray*}
(id \otimes \Delta_\sharp) \Delta_\sharp(ab)&:=& a_{(1)} \otimes \Delta_\sharp(a_{(2)}b) + ab_{(1)} \otimes \Delta_\sharp(b_{(2)}) - a_{(1)} a_{(2)} \otimes \Delta_\sharp(b), \\
&:=& a_{(1)} \otimes \Delta_\sharp(a_{(2)})b +  a_{(1)} \otimes a_{(2)}b_{(1)} \otimes b_{(2)} - a_{(1)} \otimes a_{(2)(1)}a_{(2)(2)} \otimes b \\
&  & + ab_{(1)} \otimes \Delta_\sharp(b_{(2)}) - a_{(1)} a_{(2)} \otimes \Delta_\sharp(b). 
\end{eqnarray*}
On the other hand,
\begin{eqnarray*}
(\Delta_\sharp \otimes id) \Delta_\sharp(ab)&:=& \Delta_\sharp(a_{(1)}) \otimes a_{(2)}b + \Delta_\sharp(ab_{(1)}) \otimes b_{(2)} - \Delta_\sharp(a_{(1)} a_{(2)}) \otimes b, \\
&:=& \Delta_\sharp(a_{(1)}) \otimes a_{(2)}b \\
& & + a_{(1)} \otimes a_{(2)}b_{(1)} \otimes b_{(2)} + a\Delta_\sharp(b_{(1)}) \otimes b_{(2)} - a_{(1)}a_{(2)} \otimes \Delta_\sharp(b) \\
& & - \Delta_\sharp(a_{(1)})a_{(2)} \otimes b - a_{(1)}\Delta_\sharp(a_{(2)}) \otimes b + a_{(1)(1)}a_{(1)(2)}\otimes a_{(2)} \otimes b.
\end{eqnarray*}
The two equations are equal since $(id \otimes \Delta_\sharp)\Delta_\sharp(x) = (\Delta_\sharp \otimes id)\Delta_\sharp(x)$, for $x=a, \ b$. Since $\Delta_\sharp$ is supposed to be coassociative on $kS$ and that $As(S)$ is the free associative algebra generated by $S$, $\Delta_\sharp$ is coassociative on the whole $As(S)$.
\eproof
\begin{exam}{\textbf{[Dendriform-Nijenhuis bialgebras from duplications]}}
Keep notation of Proposition \ref{TDN Iff}. Define the coproduct $\Delta: kS \xrightarrow{} kS^{\otimes 2}$ as follows:
$$ \Delta(s):= s \otimes s.$$
Then, $(As(S), \ \Delta_\sharp)$ is a dendriform-Nijenhuis bialgebra.
For instance, $$\Delta(s_1s_2)= s_1 \otimes s_1s_2 + s_1s_2 \otimes s_2 -s_1^2 \otimes s_2, \ \ \forall s_1,s_2 \in S.$$
\end{exam}
\section{Dendriform-Nijenhuis algebras}
\begin{defi}{[Dendriform-Nijenhuis algebra]}
\label{DenNin}
A {\it{Dendriform-Nijenhuis algebra}} is a $k$-vector space $DN$ equipped with nine operations $\nearrow, \ \searrow, \ \swarrow, \ \nwarrow, \ \uparrow, \ \downarrow, \ \tilde{\prec}, \ \tilde{\succ}, \ \tilde{\bullet}: DN^{\otimes 2} \xrightarrow{} DN$
verifying 28 relations. To ease notation, 7 sums operations are introduced:
\begin{eqnarray*}
x \triangleleft y &:=& x \nwarrow y + x \swarrow y + x \tilde{\prec} y, \\
x \triangleright y &:=&  x \nearrow y + x \searrow y + x \tilde{\succ} y, \\
x \bar{\bullet} y &:=&  x\uparrow y + x \downarrow y + x \tilde{\bullet} y, \\ 
x \wedge y &:=&  x \nwarrow y + x \nearrow y + x \uparrow y,\\
x \vee y &:=& x \searrow y + x \swarrow y + x \downarrow y,\\ 
x \tilde{\star} y &:=&  x \tilde{\prec} y + x \tilde{\succ} y + x \tilde{\bullet} y,
\end{eqnarray*}
and
$$
x \bar{\star} y := x \nearrow y + x \searrow y + x \swarrow y + x \nwarrow y + x \uparrow y + x \downarrow y + x \tilde{\prec} y +  x \tilde{\succ} y + x \tilde{\bullet} y.$$
That is,
$$ x \bar{\star} y := x \triangleleft y + x \triangleright y + x \bar{\bullet} y := x \wedge y + x \vee y + x \tilde{\star} y,$$
for all $x,y,z \in A$.
The 28 relations are presented in two matrices. The first one is a $7 \times 3$-matrix denoted by $(M^1_{ij})_{(i:= 1, \ldots, 7; \ j:= 1, \ldots, 3)}$. The second one is a $7 \times 1$-matrix denoted by $(M^2_{i})_{(i:= 1, \ldots, 7)}$. 

$
\begin{array}{ccccccccc}
(x \nwarrow y)\nwarrow z  &=& x \nwarrow (y \bar{\star} z); & 
(x \nearrow y)\nwarrow z &=& x \nearrow (y \triangleleft z); &
(x \wedge y)\nearrow z &=& x \nearrow (y \triangleright z); \\
(x \swarrow y)\nwarrow z &=& x \swarrow (y \wedge z); &
(x \searrow y)\nwarrow z &=& x \searrow (y \nwarrow z); &
(x \vee y)\nearrow z &=& x \searrow (y \nearrow z); \\
(x \triangleleft y)\swarrow z &=& x \swarrow (y \vee z); &
(x \triangleright y)\swarrow z &=& x \searrow (y \swarrow z); &
(x \bar{\star} y)\searrow z &=& x \searrow (y \searrow z);  \\
& & \\
(x \swarrow y)\tilde{\prec} z &=& x \swarrow (y \tilde{\star} z); & 
(x \searrow y)\tilde{\prec}  z &=& x \searrow (y \tilde{\prec} z); &
(x \vee y)\tilde{\succ} z &=& x \searrow (y \tilde{\succ} z);  \\
(x \nwarrow y)\tilde{\prec} z &=& x \tilde{\prec} (y \vee z); &
(x \nearrow y)\tilde{\prec} z &=& x \tilde{\succ} (y \swarrow z); &
(x \wedge y)\tilde{\succ}  z &=& x \tilde{\succ} (y \searrow z);  \\
(x \tilde{\prec} y)\nwarrow z &=& x \tilde{\prec} (y \wedge z); &
(x \tilde{\succ} y)\nwarrow z &=& x \tilde{\succ} (y \nwarrow z); &
(x \tilde{\star} y)\nearrow z &=& x \tilde{\succ} (y \nearrow z);  \\
& & \\
(x \tilde{\prec} y) \tilde{\prec} z &=& x \tilde{\prec} (y \tilde{\star} z); &
(x \tilde{\succ} y) \tilde{\prec} z &=& x \tilde{\succ} (y \tilde{\prec} z); &
(x \tilde{\star} y) \tilde{\succ} z &=& x \tilde{\succ} (y  \tilde{\succ} z).
\end{array}
$ 
\\[0.8cm]
\begin{eqnarray*}
(x \wedge  y) \uparrow z + (x \uparrow y) \nwarrow z &=&  x \nearrow (y \bar{\bullet} z) + x \uparrow (y \bar{\star} z); \\ 
(x \vee  y) \uparrow z + (x \downarrow y)\nwarrow  z &=&  x\searrow  (y \uparrow z) + x \downarrow  (y \wedge z); \\ 
(x \bar{\star} y)\downarrow  z + (x \bar{\bullet} y)\swarrow  z &=&  x \searrow (y  \downarrow z) + x \downarrow (y \vee z); \\ 
& & \\
(x \vee y) \tilde{\bullet} z + (x \downarrow  y)\tilde{\prec}  z &=&  x \searrow (y \tilde{\bullet} z) + x \downarrow  (y \tilde{\star} z); \\ 
(x \wedge y) \tilde{\bullet} z + (x \uparrow y)\tilde{\prec}  z &=&  x\tilde{\succ}  (y \downarrow z) + x \tilde{\bullet} (y \vee z); \\ 
(x \tilde{\star} y) \uparrow z + (x \tilde{\bullet} y)\nwarrow  z &=&  x \tilde{\succ} (y \uparrow z) + x \tilde{\bullet} (y \wedge z); \\
& & \\
(x \tilde{\star} y) \tilde{\bullet} z + (x \tilde{\bullet} y)\tilde{\prec}  z &=&  x \tilde{\succ} (y \tilde{\bullet} z) + x \tilde{\bullet} (y \tilde{\star} z). 
\end{eqnarray*} 
The {\it{vertical structure}} of the dendriform-Nijenhuis algebra $DN$ is by definition the $k$-vector space $DN_v := (DN, \ \triangleleft, \ \triangleright, \ \bar{\bullet})$.  Its {\it{horizontal structure}}  is by definition the $k$-vector space $DN_h := (DN, \ \wedge, \ \vee , \ \tilde{\star})$.
\end{defi}
\Rk \textbf{[Opposite structure of a dendriform-Nijenhuis algebra]}
There exists a symmetry letting the two matrices of relations of a dendriform-Nijenhuis algebra globally invariant. This allows the construction of the so called {\it{opposite structure}}, defined as follows.
$$x \searrow^{op} y = y \nwarrow x, \ \ x \nearrow^{op} y = y \swarrow x, \ \ x \nwarrow^{op} y = y \searrow x, \ \
x \swarrow^{op} y = y \nearrow x, \ \ x \uparrow^{op} y = y \downarrow x,$$
$$x \downarrow^{op} y = y \uparrow x, \ \ x \tilde{\succ}^{op} y =  y \tilde{\prec} x, \ \ x \tilde{\bullet}^{op} y =y \tilde{\bullet} x, \ \ x \tilde{\prec}^{op} y =  y \tilde{\succ} x .$$
Therefore, $x \triangleleft^{op} y = y \triangleright x, \ \ x \triangleright^{op} y = y \triangleleft x, \ \ x \vee^{op} y = y \wedge x, \ \ x \wedge^{op} y = y \vee x$.
A dendriform-Nijenhuis algebra is said to be {\it{commutative}} when it coincides with its opposite. For any $x,y \in T$, observe that
$x \bar{\star} y := y \bar{\star} x$. Observe that the matrix of relations $M^1$ has 3 centers of symmetry. The first one, $M^1 _{22}$,
corresponds to the first bloc of three rows, the second one, $M^1 _{52}$, to the second bloc of three rows. The last one is $M^1 _{72}$.
There are also 3 centers of symmetry for the matrix of relations $M^2$.
The first one, $M^2 _{2}$,
corresponds to the first bloc of three rows, the second one, $M^2 _{5}$,
to the second bloc of three rows. The last one is $M^2 _{7}$.
\begin{theo}
\label{verhon}
Keep notation of Definition \ref{DenNin}. Let $DN$ be a dendriform-Nijenhuis algebra. Then, its vertical structure $DN_v := (DN, \ \triangleleft, \ \triangleright, \ \bar{\bullet})$ is a $NS$-algebra
and its horizontal structure $DN_h := (DN, \ \wedge, \ \vee , \ \bar{\star})$ is a dendriform trialgebra. 
\end{theo}
\Proof
Let $DN$ be a dendriform-Nijenhuis algebra. The
vertical structure $DN_v := (DN, \ \triangleleft, \ \triangleright, \ \bar{\bullet})$ is a $NS$-algebra. Indeed for all $x,y,z \in DN$,
\begin{eqnarray*}
\sum_{i=1,2,3,7} M^1 _{i1} - \sum_{i=4,5,6} M^1 _{i1} &\Leftrightarrow &  (x \triangleleft y) \triangleleft z = x \triangleleft (y \bar{\star} z), \\
\sum_{i=1,2,3,7} M^1 _{i2} - \sum_{i=4,5,6} M^1 _{i2}  & \Leftrightarrow & (x \triangleright y) \triangleleft z = x \triangleright (y \triangleleft z),\\
\sum_{i=1,2,3,7} M^1 _{i3} - \sum_{i=4,5,6} M^1 _{i3} & \Leftrightarrow & (x \bar{\star} y) \triangleright z = x \triangleright (y \triangleright z),\\
\sum_{i=1,2,3,7} M^2 _{i} - \sum_{i=4,5,6} M^2 _{i} & \Leftrightarrow & (x \bar{\star} y) \bar{\bullet} z + (x \bar{\bullet} y) \triangleleft z = x \triangleright (y \bar{\bullet} z) + x \bar{\bullet} (y \bar{\star} z).
\end{eqnarray*}
The
horizontal structure $DN_h := (DN, \ \triangleleft, \ \triangleright, \ \bar{\bullet})$ is a dendriform trialgebra. Indeed for all $x,y,z \in DN$,
\begin{eqnarray*}
\sum_{j=1,2,3} M^1 _{1j} - M^2 _{1} &\Leftrightarrow &  (x \wedge y) \wedge z = x \wedge (y \bar{\star} z), \\
\sum_{j=1,2,3,} M^1 _{2j} -  M^2 _{2}  & \Leftrightarrow & (x \vee y) \wedge z = x \vee (y \wedge z),\\
\sum_{j=1,2,3} M^1 _{3j} -  M^2 _{3} & \Leftrightarrow & (x \bar{\star} y) \vee z = x \vee (y \vee z),\\
\sum_{j=1,2,3} M^1 _{4j} -  M^2 _{4} & \Leftrightarrow & (x \vee y) \tilde{\star}  z = x \vee (y \tilde{\star} z),\\
\sum_{j=1,2,3} M^1 _{5j} -  M^2 _{5} & \Leftrightarrow & (x \wedge  y)  \tilde{\star} z = x \tilde{\star} (y \vee z),\\
\sum_{j=1,2,3} M^1 _{6j} -  M^2 _{6} & \Leftrightarrow & (x \tilde{\star} y) \wedge z = x \tilde{\star} (y \wedge z),\\
\sum_{j=1,2,3} M^1 _{7j} -  M^2 _{7} & \Leftrightarrow & (x \tilde{\star} y) \tilde{\star} z = x \tilde{\star}  (y \tilde{\star} z).\\
\end{eqnarray*}
\eproof
\Rk
In a categorical point of view, Theorem \ref{verhon} gives two functors $F_v$ and $F_h$ represented in the following diagram:
\begin{center}

$
\begin{array}{ccc}
\textsf{Dend. Nijenhuis} & \stackrel{ F_v}{ \longrightarrow} & \textsf{NS} \\ 
 & & \\
 F_h  \downarrow & \searrow   f  & \downarrow  F_1  \\
& & \\
\textsf{TriDend.} & \stackrel{F_2}{\longrightarrow}  & \textsf{As.} 
\end{array} 
$
\end{center}
This diagram commutes, i.e., $F_2F_h = f= F_1F_v$.
\begin{defi}{}
\label{Ninop}
Let $(\mathcal{N}, \ \prec, \ \succ, \ \bullet)$ be a $NS$-algebra. 
Set $N(2) :=  k \{\ \prec, \ \succ, \ \bullet \ \}$.
A $TD$-operator $\gamma$ on $\mathcal{N}$ is a linear map $\gamma: \mathcal{N} \xrightarrow{} \mathcal{N}$ such that, 
\begin{enumerate}
\item {There exit $i \in \mathcal{N}$ and  
two linear maps,
$$
\centerdot_1: \mathcal{N} \otimes k \bra \gamma(i) \ket   \xrightarrow{} \mathcal{N}, \ \  \textrm{and } \ \  \centerdot_2: k \bra \gamma(i) \ket \otimes \mathcal{N} \xrightarrow{} \mathcal{N}
,$$
$$ 
x  \otimes \lambda \gamma(i) \mapsto \lambda x  \centerdot_1 \gamma(i),
\ \ \textrm{and } \ \  \lambda \gamma(i) \otimes x \mapsto \lambda \gamma(i) \centerdot_2 x,$$
such that
$x \diamond \gamma(i) \centerdot_2 y = x \centerdot_1 \gamma(i) \diamond y := -x \tilde{\diamond} y $, for all $\diamond \in N(2)$ and $x,y \in \mathcal{N}$,}
\item {In addition, $$\gamma(x) \diamond \gamma(y) = \gamma(\gamma(x) \diamond y + x \diamond \gamma(y) +x \tilde{\diamond}  y), $$
for all $x,y \in \mathcal{N}$ and $\gamma(x) \centerdot_1 \gamma(i) = \gamma(i) \centerdot_2 \gamma(x)$.}
\item {For all $\diamond \in N(2), \ x,y \in \mathcal{N}$, $\gamma(i) \centerdot_2 \gamma(x) \diamond y =\gamma(i) \centerdot_2 (\gamma(x) \diamond y)$ and $ x \diamond \gamma(y)\centerdot_1 \gamma(i) = (x \diamond \gamma(y))\centerdot_1 \gamma(i).$}
\end{enumerate}
\end{defi}
\begin{prop}
\label{Ninprop}
Let $(A, \ \mu)$ be a unital associative algebra with unit $i$. Suppose there exists a Nijenhuis operator $\beta: (A, \ \mu) \xrightarrow{} (A, \ \mu) $ which commutes with a $TD$-operator $\gamma: (A, \ \mu) \xrightarrow{} (A, \ \mu) $. Suppose $\gamma(i) = \beta(i)$. Then,
$\gamma$ is a $TD$-operator on the $NS$-algebra $A^{\beta}$.
\end{prop}
\Proof
By proposition \ref{asnij}, $A^{\beta}$ is a $NS$-algebra. By applying Proposition \ref{gamnin}, items 1 and 3 of Definition \ref{Ninop} hold since $\gamma(i) = \beta(i)$ and $\gamma(i)\gamma(x)= \gamma(x)\gamma(i)$, for all $x \in A$.
Fix $x,y \in A$. 
\begin{eqnarray*}
\gamma(x) \prec_{\beta} \gamma(y) &=& \gamma(x) \beta( \gamma(y)) = \gamma(x)  \gamma(\beta(y)), \\
&=& \gamma(\gamma(x)  \beta(y) + x  \gamma(\beta(y)) -x  \gamma(i) \beta(y)), \\
&=& \gamma(\gamma(x) \prec_{\beta} y + x \prec_{\beta} \gamma(y) +x \tilde{\prec}_{\beta} y). 
\end{eqnarray*}
Checking the three other equations is straightforward.
\eproof
\Rk 
Observe that $\gamma(x) \star_{\beta} \gamma(y) = \gamma(\gamma(x) \star_{\beta} y + x \star_{\beta} \gamma(y) +x \tilde{\star}_{\beta}  y)$ and thus $x \bar{\star} y := \gamma(x) \star_{\beta} y + x \star_{\beta} \gamma(y) +x \tilde{\star}_{\beta}  y$ is an associative product, or that $\gamma: (A, \ \bar{\star}) \xrightarrow{} (A, \ \star_{\beta})$ is a morphism of associative algebras.
\begin{prop}
\label{trans}
Let $(\mathcal{N}, \ \prec, \ \succ, \ \bullet)$ be a $NS$-algebra, with
$\star \longrightarrow  \ \prec + \succ + \bullet$ and $\gamma$ be a $TD$-operator on $\mathcal{N}$. Denote by $i \in \mathcal{N}$ the element which verifies items 1,2 and 3 of Definition \ref{Ninop}.
For all $x,y \in \mathcal{N}$, define nine operations as follows,
$$ x \searrow_\gamma y = \gamma(x) \succ y, \ \
x \nwarrow_\gamma y = x \prec \gamma(y), \ \
x \nearrow_\gamma y = x \succ \gamma(y), \ \
x \swarrow_\gamma y = \gamma(x) \prec y, \ \ $$
$$ x \uparrow_\gamma y = x \bullet \gamma(y), \ \ 
 x \downarrow_\gamma y = \gamma(x) \bullet y,$$
and,
$$ x \tilde{\prec }_\gamma y :=- x \gamma(i)\prec y = -x \prec \gamma(i)y, \ \  x \tilde{\succ}_\gamma  y :=  -x \succ \gamma(i) y =- x \gamma(i) \succ  y, \ \  x \tilde{\bullet}_\gamma y := -x \gamma(i)\bullet y =- x \bullet \gamma(i) y, $$ for all $x,y \in A.$ 
Define also sums operations as follows,
\begin{eqnarray*}
 x \triangleright_{\gamma} y &:=& x \searrow_{\gamma} y + x \nearrow_{\gamma} y +  x \tilde{\succ}_{\gamma} y := \gamma(x) \succ y + x \succ \gamma(y) + x \tilde{\succ}_{\gamma} y, \\
 x \triangleleft_{\gamma} y &:=& x \nwarrow_{\gamma} y + x \swarrow_{\gamma} y + x \tilde{\prec}_{\gamma} y := \gamma(x) \prec y + x \prec \gamma(y) + x \tilde{\prec}_{\gamma} y, \\
x \bar{\bullet}_{\gamma} y &:=& x \uparrow_{\gamma} y + x \downarrow_{\gamma} y + x \tilde{\bullet}_{\gamma} y :=  x \bullet \gamma(y) + \gamma(x) \bullet y + x \tilde{\bullet}_{\gamma} y, \\
x \vee_{\gamma} y &:=& x \searrow_{\gamma} y + x \swarrow_{\gamma} y +  x \downarrow_{\gamma} y := \gamma(x) \star y, \\
x \wedge_{\gamma} y &:=& x \nearrow_{\gamma} y + x \nwarrow_{\gamma} y +  x \uparrow_{\gamma} y:= x \star \gamma(y),\\
x \tilde{\star}_{\gamma} y &:=&  x \tilde{\prec}_{\gamma} y +  x \tilde{\succ}_{\gamma} y  +  x \tilde{\bullet}_{\gamma} y, \\
x \bar{\star}_{\gamma} y &:=& x \vee_{\gamma} y + x \wedge_{\gamma} y + x \tilde{\star}_{\gamma} y = x \triangleright_{\gamma} y + x \triangleleft_{\gamma} y + x \bar{\bullet}_{\gamma} y.
\end{eqnarray*}
Equipped with these nine operations, $\mathcal{N}$ is a dendriform-Nijenhuis algebra.
\end{prop}
\Proof
Let us check the relation $M^1_{11}:= (x \nwarrow_\gamma y)\nwarrow_\gamma z  = x \nwarrow_\gamma (y \bar{\star}_\gamma z)$ of Definition \ref{DenNin}.
\begin{eqnarray*}
(x \nwarrow_\gamma y)\nwarrow_\gamma z  &=& (x \prec \gamma(y)) \prec \gamma(z); \\
&=& x \prec (\gamma(y) \star \gamma(z)); \\
&=& x \prec \gamma(y \bar{\star}_\gamma z); \\
&=& x \nwarrow_\gamma (y \bar{\star}_\gamma z).
\end{eqnarray*}
Similarly, let us check the relation $M^1_{51} := 
(x \nwarrow_\gamma y)\tilde{\prec}_\gamma z = x \tilde{\prec}_\gamma (y \vee_\gamma z)$.
\begin{eqnarray*}
(x \nwarrow_\gamma y)\tilde{\prec}_\gamma z &=& (x \prec \gamma(y)) \prec \gamma(i)\centerdot_2 z,\\
&=& x \prec (\gamma(y) \star \gamma(i) \centerdot_2 z), \\
&=& x \prec (\gamma(y) \centerdot_1 \gamma(i) \star z), \\
&=& x \prec (\gamma(i) \centerdot_2 \gamma(y) \star z), \\
&=& x \prec \gamma(i) \centerdot_2 ( \gamma(y) \star z), \\
&=& x \tilde{\prec}_\gamma (y \vee_\gamma z), \\
\end{eqnarray*}
Checking the 26 other axioms does not present any difficulties. 
\eproof
\section{Transpose of a dendriform-Nijenhuis algebra}
\begin{defi}{[Transpose of a dendriform-Nijenhuis algebra]}
\label{deftrans}
A dendriform-Nijenhuis algebra $DN_1:=(DN, \searrow_1, \ \nwarrow_1,\ \nearrow_1, \ \swarrow_1, \ \uparrow_1, \ \downarrow_1, \ \tilde{\prec}_1, \ \tilde{\succ}_1, \ \tilde{\bullet}_1)$ is said to be the {\it{transpose}} of a dendriform-Nijenhuis algebra $DN_2:=(DN, \searrow_2, \ \nwarrow_2,\ \nearrow_2, \ \swarrow_2, \ \uparrow_2, \ \downarrow_2, \ \tilde{\prec}_2, \ \tilde{\succ}_2, \ \tilde{\bullet}_2)$ if for all $x,y \in DN$,
$$
x \searrow_1 y = x \searrow_2 y, \ \ 
x \nwarrow_1 y = x \nwarrow_2 y, \ \ 
x \nearrow_1 y = x \swarrow_2 y, \ \ 
x \swarrow_1 y = x \nearrow_2 y, \ \
x \uparrow_1 y = x \tilde{\prec}_2 y, \ \
$$
$$
x \downarrow_1 y = x \tilde{\succ}_2 y, \ \
x \tilde{\prec}_1 y = x \uparrow_2 y, \ \
x \tilde{\succ}_1 y = x \downarrow_2 y, \ \
x \tilde{\bullet}_1 y = x \tilde{\bullet}_2 y, \ \
x \triangleright_1 y = x \vee_2 y, \ \
$$
$$
x \triangleleft_1 y = x \wedge_2 y, \ \
x \vee_1 y = x \triangleright_2 y, \ \
x \wedge_1 y = x \triangleleft_2 y, \ \
x \bar{\bullet}_1 y = x \tilde{\star}_2 y, \ \
x \tilde{\star}_1 y = x \bar{\bullet}_2 y.
$$ 
\end{defi}
\begin{defi}{}
A {\it{Nijenhuis operator on a dendriform trialgebra}} $(TD, \ \prec, \ \succ, \ \circ)$ is a linear map $\beta: TD \xrightarrow{} TD$ such that
for all $x,y \in TD$ and $\diamond \in \{\prec, \ \succ, \ \circ \}$,
$$ \beta(x) \diamond \beta(y) = \beta(\beta(x)\diamond  y + x \diamond \beta(y) - \beta(x \diamond y)).$$
\end{defi}
\Rk
If $\star$ denotes the associative operation of a dendriform trialgebra, then $\beta(x) \star \beta(y) = \beta(\beta(x)\star  y + x \star \beta(y) - \beta(x \star y))$. Therefore, the map $\beta$ is a Nijenhuis operator on the associative algebra $(TD, \ \star)$ and a morphism of associative algebra $(A, \ \bar{\star}) \xrightarrow{} (A, \ \star)$, where the associative operation $\bar{\star}$ is defined by $x \bar{\star} y:= \beta(x)\star  y + x \star \beta(y) - \beta(x \star y)$, for all $x,y \in TD$.
\begin{prop}
\label{Nijtrans}
Let $A$ be a unital associative algebra. Suppose $\beta: A \xrightarrow{} A$
is a Nijenhuis operator which commutes with a $TD$-operator $\gamma: A \xrightarrow{} A$. Then, $\beta$
is a Nijenhuis operator on the dendriform trialgebra $A^{\gamma}$.
\end{prop}
\Proof
By Proposition \ref{TDdend}, $A^{\gamma}$ is a dendriform trialgebra. Fix 
$x,y \in A$. For instance,
\begin{eqnarray*}
\beta(x) \prec_{\gamma} \beta(y)&:=& \beta(x) \gamma( \beta(y)), \\
&=& \beta(x)  \beta(\gamma(y)), \\
&=& \beta(\beta(x)\gamma(y) + x  \beta(\gamma(y)) - \beta(x \gamma(y))), \\
&=& \beta(\beta(x)\gamma(y) + x  \gamma(\beta(y)) - \beta(x \gamma(y))), \\
&=& \beta(\beta(x) \prec_{\gamma} y + x  \prec_{\gamma} \beta(y) - \beta(x \prec_{\gamma} y)). 
\end{eqnarray*}
\eproof
\begin{prop}
\label{trans1}
Let $(TD, \ \prec, \ \succ, \ \circ)$ be a dendriform trialgebra and $\beta$ be a Nijenhuis operator. 
For all $x,y \in TD$, define nine operations as follows,
$$ x \searrow_\beta y = \beta(x) \succ y, \ \
x \nwarrow_\beta y = x \prec \beta(y), \ \
x \nearrow_\beta y = x \succ \beta(y), \ \
x \swarrow_\beta y = \beta(x) \prec y, \ \ $$
$$ x \uparrow_\beta y = x \circ \beta(y), \ \ 
 x \downarrow_\beta y = \beta(x) \circ y,$$
and,
$$ x \tilde{\prec}_\beta y := -\beta(x \prec y), \ \  x \tilde{\succ}_\beta y :=  -\beta(x \succ  y), \ \  x \tilde{\bullet}_\beta y := -\beta(x \circ y), \ \ \forall \ x,y \in A.$$ 
The sum operations are defined as in Definition \ref{DenNin}.
Then, for all $x,y,z \in TD$, (the label $\beta$ and  being omitted),

$
\begin{array}{cccccccccc}
(x \nwarrow y)\nwarrow z &=& x \nwarrow (y \bar{\star} z),&
 (x \swarrow y)\nwarrow z &=& x \swarrow (y \wedge z),&
(x \triangleleft y)\swarrow z &=& x \swarrow (y \vee z),  \\
(x \nearrow y)\nwarrow z &=& x \nearrow (y \triangleleft z),&
(x \searrow y)\nwarrow z &=& x \searrow (y \nwarrow z),&
(x \triangleright y)\swarrow z &=& x \searrow (y \swarrow z), \\
(x \wedge y)\nearrow z &=& x \nearrow (y \triangleright z),&
(x \vee y)\nearrow z &=& x \searrow (y \nearrow z),&
(x \bar{\star} y)\searrow z &=& x \searrow (y \searrow z), \\ 
& &\\
(x \nearrow y)\uparrow z &=& x \nearrow (y \bar{\bullet} z), &
(x \searrow y)\uparrow z &=& x \searrow (y \uparrow z),&
(x \triangleright y)\downarrow z &=& x \searrow (y \downarrow z), \\
(x \nwarrow y)\uparrow z &=& x \uparrow (y \triangleright z),&
(x \swarrow y)\uparrow z &=& x \downarrow (y \nearrow z),& 
(x \triangleleft y)\downarrow z &=& x \downarrow (y \searrow z), \\
(x \uparrow y)\nwarrow z &=& x \uparrow (y \triangleleft z),&
(x \downarrow y)\nwarrow z &=& x \downarrow (y \nwarrow z),&
(x \bar{\bullet} y)\swarrow z &=& x \downarrow (y \swarrow z),  \\
(x \uparrow y)\uparrow z &=& x \uparrow (y \bar{\bullet} z),&
(x \downarrow y)\uparrow z &=& x \downarrow (y \uparrow z),&
(x \bar{\bullet} y)\downarrow z &=& x \downarrow (y \downarrow z).
\end{array} 
$

\begin{eqnarray*}
(x \triangleleft y)\tilde{\prec}  z + (x \tilde{\prec} y) \nwarrow z &=& x \swarrow (y \tilde{\star} z) + x \tilde{\prec} (y \bar{\star} z), \\
(x \triangleright y) \tilde{\prec}  z + (x \tilde{\succ} y) \nwarrow z &=& x \searrow (y \tilde{\prec} z) + x \tilde{\succ} (y \triangleleft z), \\
(x \bar{\star} y) \tilde{\succ} z + (x \tilde{\star} y) \nearrow z &=& x \searrow (y \tilde{\succ} z) + x \tilde{\succ} (y \triangleright z), \\
(x \triangleright y) \tilde{\bullet}  z + (x \tilde{\succ} y) \uparrow  z &=& x \searrow (y \tilde{\bullet} z) + x  \tilde{\succ} (y \bar{\bullet} z), \\
(x \triangleleft y) \tilde{\bullet} z + (x \tilde{\prec} y) \uparrow z &=& x \downarrow (y \tilde{\succ} z) + x \tilde{\bullet} (y  \triangleright z), \\
(x \bar{\bullet} y) \tilde{\prec} z + (x \tilde{\bullet} y) \nwarrow z &=& x \downarrow (y \tilde{\prec} z) + x \tilde{\bullet} (y \triangleleft z), \\
(x \bar{\bullet} y) \tilde{\bullet} z + (x \tilde{\bullet} y) \uparrow z &=& x \downarrow (y \tilde{\bullet} z) + x  \tilde{\bullet} (y \bar{\bullet} z).
\end{eqnarray*}
Otherwise stated, $(TD, \searrow_{\beta}, \ \nwarrow_{\beta},\ \nearrow_{\beta}, \ \swarrow_{\beta}, \ \uparrow_{\beta}, \ \downarrow_{\beta}, \ \tilde{\prec}_{\beta}, \ \tilde{\succ}_{\beta}, \ \tilde{\bullet}_{\beta})$ is a dendriform-Nijenhuis algebra where 
$\searrow_{\beta}$ plays the r\^ole of $\searrow$,
$\nwarrow_{\beta} \equiv \nwarrow$,  
$ \nearrow_{\beta}  \equiv  \swarrow$, 
$ \swarrow_{\beta}  \equiv  \nearrow$, 
$ \uparrow_{\beta}  \equiv \tilde{\prec}$, 
$ \downarrow_{\beta}  \equiv  \tilde{\succ}$,
$ \tilde{\prec}_{\beta} \equiv  \uparrow$,
$ \tilde{\succ}_{\beta} \equiv  \downarrow$,
$ \tilde{\bullet}_{\beta} \equiv  \tilde{\bullet}$,
$ \triangleright_{\beta} \equiv  \vee$,
$ \triangleleft_{\beta}  \equiv \wedge$,
$ \vee_{\beta} \equiv  \triangleright$,
$ \wedge_{\beta} \equiv \triangleleft$,
$ \bar{\bullet}_{\beta} \equiv  \tilde{\star}$, and
$ \tilde{\star}_{\beta}  \equiv \bar{\bullet}$.
\end{prop}
\Proof
Keep notation of Proposition \ref{trans1}.
Fix $x,y,z \in TD$. 
Let us check for instance, $$(x \uparrow y)\uparrow z = x \uparrow (y \bar{\bullet} z),$$ the label $\beta$ being omitted on the operations.
\begin{eqnarray*}
(x \uparrow y)\uparrow z &=& (x \circ \gamma(y))\circ \gamma(z), \\
&=& x \circ (\gamma(y) \circ \gamma(z)), \\
&=& x \circ \gamma(y \bar{\bullet} z), \\
&=& x \uparrow (y \bar{\bullet} z). \\
\end{eqnarray*}
\eproof
\begin{prop}
\label{commut}
Let $(A, \mu)$ be a unital associative algebra with unit $i$. Suppose
$\beta: A \xrightarrow{} A$ is a Nijenhuis operator which commutes with
a $TD$-operator $\gamma: A \xrightarrow{} A$ and $\gamma(i)=\beta(i)$. Then, the dendriform-Nijenhuis algebra
obtained by action of $\beta$ on the dendriform trialgebra $A^{\gamma}$ is the transpose of the dendriform-Nijenhuis algebra
obtained by action of $\gamma$ on the $NS$-algebra $A^{\beta}$.
\end{prop}
\Proof
By Proposition \ref{TDdend}, the action of the $TD$-operator $\gamma$ on $A$ yields a dendriform trialgebra $A^{\gamma}$. By Proposition \ref{Nijtrans}, $\beta$ which commutes with $\gamma$ is a Nijenhuis operator on $A^{\gamma}$. By Proposition \ref{trans1}, $A$ has a dendriform-Nijenhuis algebra structure. Conversely, By Proposition \ref{asnij}, the action of the Nijenhuis operator $\beta$ on $A$ yields a Nijenhuis operator $A^{\beta}$. By Proposition \ref{Ninprop}, $\gamma$ is a $TD$-operator on $A^{\beta}$. The $k$-vector space $A$ has then another 
dendriform-Nijenhuis algebra structure by Proposition \ref{trans}.
Fix $x,y \in A$ and keep notation of Propositions \ref{trans} and \ref{trans1}. Observe that, $x \nearrow_{\gamma} y := x \succ_{\beta} \gamma(y) := \beta(x)\gamma(y)$ and $x \swarrow_{\beta} y:= \beta(x) \prec_{\gamma} y := \beta(x)\gamma(y) $. Therefore, $x \nearrow_{\gamma} y = x \swarrow_{\beta} y$ as expected in Definition \ref{deftrans}.
Similarly, $x \downarrow_{\gamma} y := \gamma(x) \bullet_{\beta} y := -\beta( \gamma(x)y )$ and $x \tilde{\succ}_{\beta} y := -\beta (x \succ_{\gamma} y) := 
-\beta( \gamma(x)y )$. Therefore, $x \downarrow_{\gamma} y:= x \tilde{\succ}_{\beta} y$. Similarly, $x \tilde{\prec}_{\gamma} y := -x \gamma(i) \prec_{\beta} y := -x \gamma(i) \beta(y) := x \circ_{\gamma} \beta(y) := x \uparrow_{\beta} y$, and so forth.
\eproof

\section{Dendriform-Nijenhuis algebras from  dendriform-Nijenhuis bialgebras}
\begin{theo}
Let $(A, \ \mu, \ \Delta)$ be a dendriform-Nijenhuis bialgebra and consider
$\textsf{End}(A)$ as an associative algebra under composition, equipped with the convolution product $*$. Then, there exists two dendriform-Nijenhuis algebra structures on $\textsf{End}(A)$, the one being the transpose of the other.
\end{theo}
\Proof
By Proposition \ref{NTDop}, the right shift $\beta$ is a Nijenhuis operator which commutes with the left shift $\gamma$ which is a $TD$-operator. If $id: A \xrightarrow{} A$ denotes the identity map, then observe that $\beta(id):= id*id := \gamma(id)$.
Use Proposition \ref{commut} to conclude.
\eproof
\Rk
Let $T,S \in \textsf{End}(A)$. By applying Proposition \ref{trans}, the nine operations are given by,
\begin{eqnarray*}
T \searrow S &:=& \beta\gamma(T)S:= (id*T*id)S, \\
T \nwarrow S &:=& T \beta\gamma(S):= T(id*S*id), \\
T \nearrow S &:=& \beta(T) \gamma(S):= (id*T)(S *id), \\
T \swarrow  S &:=& \gamma(T) \beta(S) := (T*id)(id*S), \\
T \uparrow  S &:=& -\beta(T \gamma(S)) := -id*(T(S *id)) ,\\
T \downarrow  S &:=& -\beta(\gamma(T) S):= -id*((T*id)S)        ,\\
T \tilde{\prec}  S &:=& -T \gamma(id)\beta(S) := -T(id*id)(id*S)  ,\\
T \tilde{\succ}  S &:=& -\beta(T)\gamma(id)S:=-(id*T)(id*id)S        ,\\
T \tilde{\bullet}  S &:=& \beta(T \gamma(id)  S):= id * (T(id*id)S).        
\end{eqnarray*}
The horizontal structure is a dendriform trialgebra given by,
\begin{eqnarray*}
T \wedge S &:=& T \nearrow S + T \nwarrow S + T \uparrow S:=(id*T)(S *id) +T(id*S*id) -  id*(T(S *id)),\\
T \vee S &:=& T \searrow S + T \swarrow S +   T \downarrow S:= (id*T*id)S + (T*id)(id*S) - id*((T*id)S), \\
T \tilde{\star} S &:=& T \tilde{\prec}  S + T \tilde{\succ}  S + T \tilde{\bullet}  S := -T(id*id)(id*S) -(id*T)(id*id)S  +   id * (T(id*id)S).        
\end{eqnarray*}
The vertical structure is a $NS$-algebra given by,
\begin{eqnarray*}
T \triangleright S &:=& T \nearrow S + T \searrow S + T \tilde{\succ} S:=
(id*T)(S *id) + (id*T*id)S -(id*T)(id*id)S,   \\
T \triangleleft S &:=& T \nwarrow S + T \swarrow S + T \tilde{\prec} S:= T(id*S*id) + (T*id)(id*S) - T(id*id)(id*S),
\\
T \bar{\bullet} S &:=& T \uparrow S + T \downarrow S + T \tilde{\bullet} S:= -id*(T(S *id)) - id*((T*id)S) + id * (T(id*id)S).
\end{eqnarray*}
The associative operation $\bar{\star}$ is the sum of the nine operations, i.e.,
\begin{eqnarray*}
T \bar{\star} S &:=& (id*T)(S *id) + (id*T*id)S -(id*T)(id*id)S \\ 
& &+ T(id*S*id) + (T*id)(id*S) - T(id*id)(id*S)\\
& &- id*(T(S *id) - id*((T*id)S) + id * (T(id*id)S).
\end{eqnarray*}
\section{Free $NS$-algebra and free dendriform-Nijenhuis algebra}
Let us recall what an operad is, see \cite{Lodayscd, Fresse, GK} for instance. 

\noindent
Let $P$ be a type of algebras, for instance the $NS$-algebras, and $P(V)$ be the free $P$-algebra on the $k$-vector space $V$. Suppose $P(V) := \oplus_{n \geq 1} P(n) \otimes _{S_n} V^{\otimes n},$ where $P(n)$ are right $S_n$-modules. Consider $P$ as an endofunctor on the category of $k$-vector spaces. The structure of the free $P$-algebra of $P(V)$ induces a natural transformation $\pi: P \tilde{\circ} P \xrightarrow{} P$ as well as $u: Id \xrightarrow{} P$ verifying usual associativity and unitarity axioms. An algebraic operad is then a triple $(P, \ \pi, \ u)$. A $P$-algebra is then a $k$-vector space $V$ together with a linear map $\pi_A: P(A) \xrightarrow{} A$ such that $\pi_A \tilde{\circ} \pi(A) = \pi_A \tilde{\circ} P(\pi_A)$ and $\pi_A \tilde{\circ} u(A) = Id_A.$
The $k$-vector space $P(n)$ is the space of $n$-ary operations for $P$-algebras. We will always suppose there is, up to homotheties, a unique $1$-ary operation, the identity, i.e.,  $P(1):= kId$ and that all possible operations are generated by composition from $P(2)$. The operad is said to be {\it{binary}}. It is said to be {\it{quadratic}} if all the relations between operations are consequences of relations described exclusively with the help of monomials with two operations. An operad is said to be {\it{non-symmetric}} if, in the relations, the variables $x,y,z$ appear in the same order. The $k$-vector space $P(n)$ can be written as  $P(n):= P'(n) \otimes k[S_n]$, where $P'(n)$ is also a $k$-vector space and $S_n$ the symmetric group on $n$ elements.
In this case, the free $P$-algebra is entirely induced by the free $P$-algebra on one generator $P(k):= \oplus_{n \geq 1} \ P'(n)$. The generating function of the operad $P$ is given by:
$$ f^{P}(x):= \sum \ (-1)^n \frac{\textrm{dim} \ P(n)}{n!} x^n := \sum \ (-1)^n \textrm{dim} \ P'(n) x^n.$$
Below, we will indicate the sequence $(\textrm{dim} \ P'(n))_{n \geq 1}$.

\subsection{On the free $NS$-algebra}
Let $V$ be a $k$-vector space. The {\it{free $NS$-algebra}} $\mathcal{N}(V)$ on $V$ is by definition,
a $NS$-algebra equipped with a map $i: \ V \mapsto \mathcal{N}(V)$ which satisfies the following universal property:
for any linear map $f: V \xrightarrow{} A$, where $A$ is a $NS$-algebra, there exists a unique $NS$-algebra morphism $\bar{f}: \mathcal{N}(V) \xrightarrow{} A$ such that $\bar{f} \circ i=f$.
The same definition holds for the {\it{free dendriform-Nijenhuis algebra}}.

Since the three operations of a $NS$-algebra have no symmetry and since compatibility axioms involve only monomials where $x, \ y$ and $z$ stay in the same order, the free $NS$-algebra is of the form:
$$ \mathcal{N}(V) := \bigoplus_{n \geq 1} \ \mathcal{N}_n \otimes V^{\otimes n} .$$
In particular, the free $NS$-algebra on one generator $x$ is $ \mathcal{N}(k) := \bigoplus_{n \geq 1} \mathcal{N}_n,$ where $\mathcal{N}_1:= kx$, $\mathcal{N}_2:= 
k(x \prec x) \oplus  k(x \succ x) \oplus k( x \bullet x).$
The space of three variables made out of three operations is of dimension $2 \times 3^2= 18$. As we have $4$ relations, the space $\mathcal{N}_3$ has a dimension equal to $18-4=14.$
Therefore, the sequence associated with the dimensions of $(\mathcal{N}_n)_{n \in \mathbb{N}}$ starts with $1, \ 3, \ 14 \ldots$
Finding the free $NS$-algebra on one generator is an open problem.  
However, we will show that the augmented free $NS$-algebra over a $k$-vector space $V$ has a connected Hopf algebra structure. Before, some preparations are needed. To be as self-contained as possible, we introduce some notation to expose a theorem due to Loday \cite{Lodayscd}. 

\noindent
Recall that a bialgebra $(H, \mu, \ \Delta, \ \eta, \ \kappa)$
is a unital associative algebra $(H, \mu, \ \eta)$ together with co-unital coassociative coalgebra $(H, \ \Delta, \ \kappa)$. Moreover, it is required that the coproduct $\Delta$
and the counit $\kappa$ are morphisms of unital algebras. A bialgebra is connected if there exists a filtration $(F_rH)_r$ such that $H = \bigcup_r F_rH$, where $F_0H := k1_H$ and for all $r$,
$$ F_rH := \{x \in H; \ \Delta(x) -1_H \otimes x - x\otimes 1_H \in F_{r-1}H \otimes F_{r-1}H \}.$$
Such a bialgebra admits an antipode. Consequently, connected bialgebras are connected Hopf algebras.  

\noindent
Let $P$ be a binary quadratic operad. By a \textit{unit action} \cite{Lodayscd}, we mean the choice of two linear applications:
$$\upsilon: P(2) \xrightarrow{} P(1), \ \ \ \ \ \varpi:P(2) \xrightarrow{} P(1),$$
giving sense, when possible, to $x \circ 1$ and $1 \circ x$, for all operations $\circ \in P(2)$ and for all $x$ in the $P$-algebra $A$, i.e.,
$x \circ 1 = \upsilon(\circ)(x)$ and $1 \circ x= \varpi(\circ)(x)$.
If $P(2)$ contains an associative operation, say $\bar{\star}$, then we require that $x \bar{\star} 1 := x := 1 \bar{\star} x$, i.e., $\upsilon(\bar{\star}) := Id := \varpi(\bar{\star})$.
We say that the unit action, or the couple $(\upsilon,\varpi)$ is {\it{compatible}} with the relations of the $P$-algebra $A$ if they still hold on $A_+:= k1 \oplus A$ as far as the terms as defined.
Let $A$, $B$ be two $P$-algebras such that $P(2)$ contains an associative operation $\bar{\star}$. Using the couple $(\upsilon,\varpi)$, we extend binary operations $\circ \in P(2)$ to the $k$-vector space $A \otimes 1.k \oplus k.1 \otimes B \oplus A \otimes B$ by requiring:
\begin{eqnarray}
(a \otimes b) \circ (a' \otimes b') &:= &(a \bar{\star} a') \otimes (b \circ b') \ \ \ \textrm{if} \ \ b \otimes b' \not= 1 \otimes 1, \\
(a \otimes 1) \circ (a' \otimes 1) &:= &(a \circ a') \otimes 1, \ \ \ \textrm{otherwise}.
\end{eqnarray}
The unit action or the couple $(\upsilon,\varpi)$ is said to be {\it{coherent}} with the relations of $P$ if $A \otimes 1.k \oplus k.1 \otimes B \oplus A \otimes B$, equipped with these operations is still a $P$-algebra. Observe that a necessary condition for having coherence is compatibility.

\noindent
One of the main interest of these two concepts is the construction of a connected Hopf algebra on the augmented free $P$-algebra. 
\begin{theo} [\textbf{Loday \cite{Lodayscd}}]
\label{Lodaych}
Let $P$ be a binary quadratic operad. Suppose there exists an associative operation in $P(2)$. Then, any unit action coherent with the relations of $P$ equips the augmented free $P$-algebra $P(V)_+$ on a $k$-vector space $V$ with a
coassociative coproduct $\Delta: P(V)_+ \xrightarrow{}  P(V)_+ \otimes P(V)_+,$
which is a $P$-algebra morphism.
Moreover,
$P(V)_+$ is a connected Hopf algebra.
\end{theo}
\Proof
See \cite{Lodayscd} for the proof. However, we reproduce it to make things clearer.
Let $V$ be a $k$-vector space and $P(V)$ be the free $P$-algebra on $V$.
Since the unit action is coherent, $P(V)_+ \otimes P(V)_+$ is a $P$-algebra. Consider the linear map $\delta: V \xrightarrow{} P(V)_+ \otimes P(V)_+$, given by $v \mapsto 1 \otimes v + v \otimes 1$. Since $P(V)$ is the free $P$-algebra on $V$, there exists a unique extension of $\delta$ to a morphism of augmented $P$-algebra $\Delta: P(V)_+ \xrightarrow{}  P(V)_+ \otimes P(V)_+$. Now, $\Delta$ is coassociative since the morphisms $(\Delta \otimes id)\Delta$ and 
$(id \otimes \Delta)\Delta$ extend the linear map $ V \xrightarrow{} P(V)_+ ^{\otimes 3}$ which maps $v$ to $1 \otimes 1 \otimes v + 1 \otimes v \otimes 1 + v \otimes 1 \otimes 1$. By unicity of the extension, the coproduct $\Delta$ is coassociative. The bialgebra we have just obtained is connected. Indeed, by definition, the free $P$-algebra $P(V)$ can be written as
$P(V) := \oplus_{n \geq 1} \ P(V)_n$, where $P(V)_n$ is the $k$-vector space of products of $n$ elements of $V$. Moreover, we have $\Delta(x \circ y ) := 1 \otimes (x \circ y )  + (x \circ y) \otimes 1 + x \otimes (1 \circ y) + y \otimes (x \circ 1)$, for all $x,y \in P(V)$ and $\circ \in P(2)$. The filtration of $P(V)_+$ is then $FrP(V)_+ = k.1 \oplus \bigoplus_{1 \leq n \leq r} P(V)_n$.
Therefore, $P(V)_+:=\cup_r \ FrP(V)_+$ and $P(V)_+$ is a connected bialgebra.
\eproof

\noindent
We will use this theorem to show that there exists a connected Hopf algebra structure on the augmented free $NS$-algebra as well as on the augmented free commutative $NS$-algebra.
\begin{prop}
Let $\mathcal{N}(V)$ be the free $NS$-algebra on a $k$-vector space $V$. Extend the binary operations $\prec, \ \succ$ and $\bullet$ to $\mathcal{N}(V)_+$ as follows:
$$ x \succ 1 := 0, \  \ 1 \succ x := x, \ \ 1 \prec x := 0, \ \ x \prec 1 := x, \ \ x \bullet 1 := 0 := 1 \bullet x,$$
So that, $x \star 1 = x = 1 \star x$ for all $x \in \mathcal{N}(V)$.
This choice is coherent.
\end{prop}
\Rk
We cannot extend the operations $\succ$ and $\prec$ to $k$, i.e., $1 \succ 1$ and $ 1 \prec 1$ are not defined.

\Proof
Keep notation introduced in that Section.
Firstly, let us show that this choice is compatible. Let $x,y,z \in \mathcal{N}(V)_+$. We have to show for instance that the relation $(x \prec y)\prec z = x \prec (y \star z)$ holds in $\mathcal{N}(V)_+$. Indeed, for $x=1$ we get $0=0$. For $y =1$ we get $ x \prec z = x \prec z$ and for $z=1$ we get $ x \prec y = x \prec y$. The same checking can be done for the 3 other equations. The augmented $NS$-algebra $\mathcal{N}(V)_+$ is then a $NS$-algebra.

\noindent
Secondly, let us show that this choice is coherent. Let $x_1,x_2, x_3, y_1, y_2, y_3 \in \mathcal{N}(V)_+$. We have to show that, for instance:
$$ ((x_1 \otimes y_1) \prec (x_2 \otimes y_2)) \prec (x_3 \otimes y_3) = (x_1 \otimes y_1) \prec ((x_2 \otimes y_2) \star (x_3 \otimes y_3)).$$
Indeed, if there exists a unique $y_i =1$, the other belonging to $\mathcal{N}(V)$, then, by definition we get:
$$ (x_1 \star x_2 \star x_3 ) \otimes ((y_1 \prec y_2) \prec  y_3) = (x_1 \star x_2 \star x_3 ) \otimes (y_1 \prec ( y_2 \star  y_3)),$$
which always holds since our choice of the unit action is compatible. If $y_1=y_2=y_3=1$, then the equality holds since our choice is compatible.
If $y_2=1$, $y_1=1$ and $y_3 \in \mathcal{N}(V)$, we get: $0=0$, similarly if $y_1=1$, $y_3=1$ and $y_2 \in \mathcal{N}(V)$. If $y_1 \in \mathcal{N}(V)$, $y_2=1$ and $y_3=1$, the two hand sides are equal to $(x_1 \bar{\star} x_2 \bar{\star} x_3 ) \otimes y_1$. Therefore, this equation holds in $\mathcal{N}(V)  \otimes 1.k \oplus k.1 \otimes \mathcal{N}(V)  \oplus \mathcal{N}(V)  \otimes \mathcal{N}(V) $. Checking the same thing with the 3 other relations shows that our choice of the unit action is coherent.
\eproof
\begin{coro}
There exists a connected Hopf algebra structure on the augmented free $NS$-algebra as well as on the augmented free commutative $NS$-algebra.
\end{coro}
\Proof
The first claim comes from the fact that our choice is coherent and from Theorem \ref{Lodaych}. For the second remark, observe that our choice is in agreement with the symmetry relations defining a commutative $NS$-algebra since for instance $x \prec^{op} 1 := 1 \succ x :=x$ and $1 \succ^{op} x := x \prec 1 :=x,$ for all $x \in \mathcal{N}(V)$.
\eproof
\subsection{On the free dendriform-Nijenhuis algebra}
The same claims hold for the free dendriform-Nijenhuis algebra. The associated operad is binary, quadratic and non-symmetric.
The free dendriform-Nijenhuis algebra on a $k$-vector spce $V$ is of the form:
$$ \mathcal{DN}(V) := \bigoplus_{n \geq 1} \ \mathcal{DN}_n \otimes V^{\otimes n} .$$
In particular, on one generator $x$, 
$$ \mathcal{DN}(k) := \bigoplus_{n \geq 1} \mathcal{N}_n,$$ where $\mathcal{N}_1:= kx$, $\mathcal{N}_2:= k(x \uparrow x) \oplus k(x\downarrow x)\oplus k(x \searrow x) \oplus k(x \nearrow x)\oplus k(x  \swarrow x) \oplus k( x \nwarrow x) \oplus 
k(x \tilde{\prec} x) \oplus  k(x \tilde{\succ} x) \oplus k( x \bullet x).$
The space of three variables made out of nine operations is of dimension $2 \times 9^2= 162$. As we have $28$ relations, the space $\mathcal{N}_3$ has a dimension equal to $162-28=134.$
Therefore, the sequence associated with the dimensions of $(\mathcal{DN}_n)_{n \in \mathbb{N}}$ starts with $1, \ 9, \ 134 \ldots$
Finding the free dendriform-Nijenhuis algebra on one generator is an open problem.  

However, there exists a connected Hopf algebra structure on the augmented free dendriform-Nijenhuis algebra as well as on the augmented free commutative dendriform-Nijenhuis algebra.
\begin{prop}
Let $\mathcal{DN}(V)$ be the free dendriform-Nijenhuis algebra on a $k$-vector space $V$. Extend the binary operations $\nwarrow$ and $\searrow$ to $\mathcal{DN}(V)_+$ as follows:
$$ x \nwarrow 1 := x, \  \ 1 \nwarrow x := 0, \ \ 1 \searrow x := x, \ \ x \searrow 1 := 0, \ \ \forall x \in \mathcal{DN}(V).$$
In addition, for any other operation, $\diamond \in \{ \nearrow, \ \swarrow, \ \uparrow, \ \downarrow, \ \tilde{\prec}, \ \tilde{\succ}, \ \tilde{\bullet}  \}$ choose: 
$ x \diamond 1 := 0 = 1 \diamond x$, for all $x \in \mathcal{DN}(V).$ Then, 
$$ x \triangleleft 1 =x, \ \ 1 \triangleright x =x, \ \ 1 \vee x := x , \ \ x \wedge 1 := x, \ \ x \bar{\star} 1 = x = 1 \bar{\star} x.$$ 
Moreover, this choice is coherent.
\end{prop}
\Rk
We cannot extend the operations $\searrow$ and $\nwarrow$ to $k$, i.e., $1 \searrow 1$ and $ 1 \nwarrow 1$ are not defined.

\Proof
Keep notation introduced in that Section.
Firstly, let us show that this choice is compatible. Let $x,y,z \in \mathcal{DN}(V)_+$. We have to show for instance that the relation $(x \nwarrow y)\nwarrow z = x \nwarrow (y \bar{\star} z)$ holds in $\mathcal{DN}(V)_+$. Indeed, for $x=1$ we get $0=0$. For $y =1$ we get $ x \nwarrow z = x \nwarrow z$ and for $z=1$ we get $ x \nwarrow y = x \nwarrow y$. We do the same thing with the 27 others and quickly found that the augmented dendriform-Nijenhuis algebra $\mathcal{DN}(V)_+$ is still a dendriform-Nijenhuis algebra.

\noindent
Secondly, let us show that this choice is coherent. Let $x_1,x_2, x_3, y_1, y_2, y_3 \in \mathcal{DN}(V)_+$. We have to show that, for instance:
$$(Eq. \ M^1_{11}) \ \ \  ((x_1 \otimes y_1) \nwarrow (x_2 \otimes y_2)) \nwarrow (x_3 \otimes y_3) = (x_1 \otimes y_1) \nwarrow ((x_2 \otimes y_2) \bar{\star} (x_3 \otimes y_3)).$$
Indeed, if there exists a unique $y_i =1$, the other belonging to $\mathcal{DN}(V)$, then, by definition we get:
$$ (x_1 \bar{\star} x_2 \bar{\star} x_3 ) \otimes (y_1 \nwarrow y_2) \nwarrow  y_3 = (x_1 \bar{\star} x_2 \bar{\star} x_3 ) \otimes y_1 \nwarrow ( y_2 \bar{\star}  y_3),$$
which always holds since our choice of the unit action is compatible. Similarly for $y_1=y_2=y_3=1$.
If $y_1=1=y_2$ and $y_3 \in \mathcal{DN}(V)$, we get: $0=0$, similarly if $y_1=1=y_3$ and $y_2 \in \mathcal{DN}(V)$. If $y_1 \in \mathcal{DN}(V)$ and $y_2=1=y_3$, the two hand sides of $(Eq. \ M^1_{11})$ are equal to $(x_1 \bar{\star} x_2 \bar{\star} x_3 ) \otimes y_1$. Therefore, $ (Eq. \ M^1_{11})$ holds in $\mathcal{DN}(V)  \otimes 1.k \oplus k.1 \otimes \mathcal{DN}(V)  \oplus \mathcal{DN}(V)  \otimes \mathcal{DN}(V) $. Checking the same thing with the 27 other relations shows that our choice of the unit action is coherent.
\eproof
\begin{coro}
There exists a connected Hopf algebra structure on the augmented free dendriform-Nijenhuis algebra as well as on the augmented free commutative dendriform-Nijenhuis algebra.
\end{coro}
\Proof
The first claim comes from the fact that our choice is coherent and from Theorem \ref{Lodaych}. For the second remark, observe that our choice is in agreement with the symmetry relations defining a commutative dendriform-Nijenhuis algebra since for instance $x \nwarrow^{op} 1 := 1 \searrow x :=x$ and $1 \searrow^{op} x := x \nwarrow 1 :=x,$ for all $x \in \mathcal{DN}(V)$.
\eproof
\section{Conclusion}
There exists another way to produce Nijenhuis operators by defining another type of bialgebras. Instead of defining the coproduct $\Delta$ by $\Delta(ab):= \Delta(a)b + a\Delta(b) - \mu(\Delta(a)) \otimes b $
on an associative algebra $(A, \ \mu)$, the following definition can be chosen: $\Delta(ab):= \Delta(a)b + a\Delta(b) - a \otimes \mu(\Delta(b)) $. The generalisation of what was written is straightforward by observing that the right shift $\beta$ becomes a $TD$-operator whereas the left shift $\gamma$ becomes a Nijenhuis operator.

\bibliographystyle{plain}
\bibliography{These}

\end{document}